\documentclass[12pt]{article}

\usepackage{amsmath,amssymb,amsthm,amscd,a4wide,upref}

\usepackage[enableskew]{youngtab}

\usepackage{hyperref}
\hypersetup{colorlinks,citecolor=blue,filecolor=black,linkcolor=blue,urlcolor=blue}
\usepackage{enumerate}
\usepackage{mathtools}
\usepackage{enumitem}

\usepackage{lmodern}     % set math font to Latin modern math
\usepackage[T1]{fontenc}
\begin{document}

% % % % % % % % % % % % % %
%\newcommand{\tma}{\textcolor{magenta}}
%\newcommand{\tre}{\textcolor{red}}
% % % % % % % % % % % % % %

\newcommand{\ad}{{\rm ad}}
\newcommand{\cri}{{\rm cri}}
\newcommand{\row}{{\rm row}}
\newcommand{\col}{{\rm col}}
\newcommand{\Ann}{{\rm{Ann}\ts}}
\newcommand{\End}{{\rm{End}\ts}}
\newcommand{\Rep}{{\rm{Rep}\ts}}
\newcommand{\Hom}{{\rm{Hom}}}
\newcommand{\Mat}{{\rm{Mat}}}
\newcommand{\ch}{{\rm{ch}\ts}}
\newcommand{\chara}{{\rm{char}\ts}}
\newcommand{\diag}{{\rm diag}}
\newcommand{\st}{{\rm st}}
\newcommand{\non}{\nonumber}
\newcommand{\wt}{\widetilde}
\newcommand{\wh}{\widehat}
\newcommand{\ol}{\overline}
\newcommand{\ot}{\otimes}
\newcommand{\la}{\lambda}
\newcommand{\La}{\Lambda}
\newcommand{\De}{\Delta}
\newcommand{\al}{\alpha}
\newcommand{\be}{\beta}
\newcommand{\ga}{\gamma}
\newcommand{\Ga}{\Gamma}
\newcommand{\ep}{\epsilon}
\newcommand{\ka}{\kappa}
\newcommand{\vk}{\varkappa}
\newcommand{\si}{\sigma}
\newcommand{\vs}{\varsigma}
\newcommand{\vp}{\varphi}
\newcommand{\ta}{\theta}
\newcommand{\de}{\delta}
\newcommand{\ze}{\zeta}
\newcommand{\om}{\omega}
\newcommand{\Om}{\Omega}
\newcommand{\ee}{\epsilon^{}}
\newcommand{\su}{s^{}}
\newcommand{\hra}{\hookrightarrow}
\newcommand{\ve}{\varepsilon}
\newcommand{\pr}{^{\tss\prime}}
\newcommand{\ts}{\,}
\newcommand{\vac}{\mathbf{1}}
\newcommand{\vacu}{|0\rangle}
\newcommand{\di}{\partial}
\newcommand{\qin}{q^{-1}}
\newcommand{\tss}{\hspace{1pt}}
\newcommand{\Sr}{ {\rm S}}
\newcommand{\U}{ {\rm U}}
\newcommand{\BL}{ {\overline L}}
\newcommand{\BE}{ {\overline E}}
\newcommand{\BP}{ {\overline P}}
\newcommand{\AAb}{\mathbb{A}\tss}
\newcommand{\CC}{\mathbb{C}\tss}
\newcommand{\KK}{\mathbb{K}\tss}
\newcommand{\QQ}{\mathbb{Q}\tss}
\newcommand{\SSb}{\mathbb{S}\tss}
\newcommand{\TT}{\mathbb{T}\tss}
\newcommand{\ZZ}{\mathbb{Z}\tss}
\newcommand{\DY}{ {\rm DY}}
\newcommand{\X}{ {\rm X}}
\newcommand{\Y}{ {\rm Y}}
\newcommand{\Z}{{\rm Z}}
\newcommand{\ZX}{{\rm ZX}}
\newcommand{\ZY}{{\rm ZY}}
\newcommand{\Ac}{\mathcal{A}}
\newcommand{\Lc}{\mathcal{L}}
\newcommand{\Mc}{\mathcal{M}}
\newcommand{\Pc}{\mathcal{P}}
\newcommand{\Qc}{\mathcal{Q}}
\newcommand{\Rc}{\mathcal{R}}
\newcommand{\Sc}{\mathcal{S}}
\newcommand{\Tc}{\mathcal{T}}
\newcommand{\Bc}{\mathcal{B}}
\newcommand{\Ec}{\mathcal{E}}
\newcommand{\Fc}{\mathcal{F}}
\newcommand{\Gc}{\mathcal{G}}
\newcommand{\Hc}{\mathcal{H}}
\newcommand{\Uc}{\mathcal{U}}
\newcommand{\Vc}{\mathcal{V}}
\newcommand{\Wc}{\mathcal{W}}
\newcommand{\Yc}{\mathcal{Y}}
\newcommand{\Cl}{\mathcal{C}l}
\newcommand{\Ar}{{\rm A}}
\newcommand{\Br}{{\rm B}}
\newcommand{\Ir}{{\rm I}}
\newcommand{\Fr}{{\rm F}}
\newcommand{\Jr}{{\rm J}}
\newcommand{\Or}{{\rm O}}
\newcommand{\GL}{{\rm GL}}
\newcommand{\Spr}{{\rm Sp}}
\newcommand{\Rr}{{\rm R}}
\newcommand{\Zr}{{\rm Z}}
\newcommand{\gl}{\mathfrak{gl}}
\newcommand{\middd}{{\rm mid}}
\newcommand{\ev}{{\rm ev}}
\newcommand{\Pf}{{\rm Pf}}
\newcommand{\Norm}{{\rm Norm\tss}}
\newcommand{\oa}{\mathfrak{o}}
\newcommand{\spa}{\mathfrak{sp}}
\newcommand{\osp}{\mathfrak{osp}}
\newcommand{\f}{\mathfrak{f}}
\newcommand{\se}{\mathfrak{s}}
\newcommand{\g}{\mathfrak{g}}
\newcommand{\h}{\mathfrak h}
\newcommand{\n}{\mathfrak n}
\newcommand{\m}{\mathfrak m}
\newcommand{\z}{\mathfrak{z}}
\newcommand{\Zgot}{\mathfrak{Z}}
\newcommand{\p}{\mathfrak{p}}
\newcommand{\sll}{\mathfrak{sl}}
\newcommand{\agot}{\mathfrak{a}}
\newcommand{\bgot}{\mathfrak{b}}
\newcommand{\qdet}{ {\rm qdet}\ts}
\newcommand{\Ber}{ {\rm Ber}\ts}
\newcommand{\HC}{ {\mathcal HC}}
\newcommand{\cdet}{{\rm cdet}}
\newcommand{\rdet}{{\rm rdet}}
\newcommand{\tr}{ {\rm tr}}
\newcommand{\gr}{ {\rm gr}\ts}
\newcommand{\str}{ {\rm str}}
\newcommand{\loc}{{\rm loc}}
\newcommand{\Gr}{{\rm G}}
\newcommand{\sgn}{ {\rm sgn}\ts}
\newcommand{\sign}{{\rm sgn}}
\newcommand{\ba}{\bar{a}}
\newcommand{\bb}{\bar{b}}
\newcommand{\bi}{\bar{\imath}}
\newcommand{\bj}{\bar{\jmath}}
\newcommand{\bk}{\bar{k}}
\newcommand{\bl}{\bar{l}}
\newcommand{\bp}{\bar{p}}
\newcommand{\hb}{\mathbf{h}}
\newcommand{\Sym}{\mathfrak S}
\newcommand{\fand}{\quad\text{and}\quad}
\newcommand{\Fand}{\qquad\text{and}\qquad}
\newcommand{\For}{\qquad\text{or}\qquad}
\newcommand{\for}{\quad\text{or}\quad}
\newcommand{\grpr}{{\rm gr}^{\tss\prime}\ts}
\newcommand{\degpr}{{\rm deg}^{\tss\prime}\tss}
\newcommand{\bideg}{{\rm bideg}\ts}

\renewcommand{\theequation}{\arabic{section}.\arabic{equation}}

\numberwithin{equation}{section}

\newtheorem{thm}{Theorem}[section]
\newtheorem{lem}[thm]{Lemma}
\newtheorem{prop}[thm]{Proposition}
\newtheorem{cor}[thm]{Corollary}
\newtheorem{conj}[thm]{Conjecture}
\newtheorem*{mthm}{Main Theorem}
\newtheorem*{mthma}{Theorem A}
\newtheorem*{mthmb}{Theorem B}
\newtheorem*{mthmc}{Theorem C}
\newtheorem*{mthmd}{Theorem D}

\theoremstyle{definition}
\newtheorem{defin}[thm]{Definition}

\theoremstyle{remark}
\newtheorem{remark}[thm]{Remark}
\newtheorem{example}[thm]{Example}
\newtheorem{examples}[thm]{Examples}

\newcommand{\bth}{\begin{thm}}
\renewcommand{\eth}{\end{thm}}
\newcommand{\bpr}{\begin{prop}}
\newcommand{\epr}{\end{prop}}
\newcommand{\ble}{\begin{lem}}
\newcommand{\ele}{\end{lem}}
\newcommand{\bco}{\begin{cor}}
\newcommand{\eco}{\end{cor}}
\newcommand{\bde}{\begin{defin}}
\newcommand{\ede}{\end{defin}}
\newcommand{\bex}{\begin{example}}
\newcommand{\eex}{\end{example}}
\newcommand{\bes}{\begin{examples}}
\newcommand{\ees}{\end{examples}}
\newcommand{\bre}{\begin{remark}}
\newcommand{\ere}{\end{remark}}
\newcommand{\bcj}{\begin{conj}}
\newcommand{\ecj}{\end{conj}}

\newcommand{\bal}{\begin{aligned}}
\newcommand{\eal}{\end{aligned}}
\newcommand{\beq}{\begin{equation}}
\newcommand{\eeq}{\end{equation}}
\newcommand{\ben}{\begin{equation*}}
\newcommand{\een}{\end{equation*}}

\newcommand{\bpf}{\begin{proof}}
\newcommand{\epf}{\end{proof}}

\def\beql#1{\begin{equation}\label{#1}}

\newcommand{\Res}{\mathop{\mathrm{Res}}}

\title{\Large\bf Representations of the super Yangians of types $A$ and $C$}

\author{Alexander Molev}

\date{} % Start October 2021
\maketitle

%\vspace{4 mm}

\begin{abstract}
We classify the finite-dimensional irreducible
representations of the super Yangian associated
with the orthosymplectic Lie superalgebra $\osp_{2|2n}$. The classification is given in terms of
the highest weights and Drinfeld polynomials.
We also include an $R$-matrix construction of the polynomial evaluation modules over the Yangian
associated with the Lie superalgebra $\gl_{m|n}$, as an appendix.
This is a super-version of the well-known
construction for the $\gl_n$ Yangian and
it relies on the Schur--Sergeev duality.

\end{abstract}

%\newpage

%\tableofcontents
%
%\newpage

\section{Introduction}\label{sec:int}
\setcounter{equation}{0}

The {\em Yangian} $\Y(\osp_{M|2n})$ associated with the orthosymplectic Lie superalgebra $\osp_{M|2n}$
is a deformation of the universal enveloping algebra $\U(\osp_{M|2n}[u])$ in the class of
Hopf algebras. The original definition is due to Arnaudon {\it et al.\/}~\cite{aacfr:rp},
where some basic properties of the Yangian were also described. Its
$\gl_{m|n}$ counterpart was introduced earlier by Nazarov~\cite{n:qb}, and
the finite-dimensional irreducible representations of the Yangian $\Y(\gl_{m|n})$
were classified by Zhang~\cite{zh:sy} in a way similar to the
Drinfeld Yangians \cite{d:nr}. Apart from the case of the Yangian $\Y(\osp_{1|2n})$
considered in \cite{m:ry}, the general
classification problem for the orthosymplectic Yangians remains open.
We show in this paper that the finite-dimensional irreducible
representations of the Yangian $\Y(\osp_{2|2n})$ are
described in essentially the same way as those of the Yangian $\Y(\gl_{m|n})$,
thus extending the similarity between the modules over the Lie superalgebras
of types $A$ and $C$ as given by Kac~\cite{k:rc}.

In more detail, we will regard the Yangian
$\Y(\osp_{2|2n})$
as a quotient
of the {\em extended Yangian} $\X(\osp_{2|2n})$ by an ideal generated by central elements.
The highest weight representation $L(\la(u))$ of $\X(\osp_{2|2n})$ is defined as
the irreducible quotient of the Verma module $M(\la(u))$
associated with an $(n+2)$-tuple $\la(u)=(\la_1(u),\dots,\la_{n+2}(u))$ of formal series
in $u^{-1}$. The tuple is called the {\em highest weight}
of the representation.

A standard argument shows that every finite-dimensional irreducible representation
of the extended Yangian is of the form $L(\la(u))$ for a certain highest weight $\la(u)$.
The key step in the classification is to determine
the necessary and sufficient conditions on the highest weight for
the representation $L(\la(u))$
to be finite-dimensional. We will rely
on reduction properties of the representations
to establish necessary conditions, and
their sufficiency will be verified by constructing a family of
{\em fundamental representations} of the Yangian $\X(\osp_{2|2n})$.
This will prove the following theorem.

\begin{mthm}\label{thm:yclassi}
Every finite-dimensional irreducible representation of the algebra $\X(\osp_{2|2n})$
is isomorphic to $L(\la(u))$ for a certain highest weight $\la(u)$.
The representation $L(\la(u))$
is finite-dimensional if and only if
there exist monic polynomials
$\ol Q(u),Q(u),P_2(u),\dots,P_{n+1}(u)$ in
$u$ such that
\begin{align}
\frac{\la_1(u)}{\la_{2}(u)}&=\frac{\ol Q(u)}{Q(u)},
\label{gloo}\\[0.4em]
\frac{\la_{i+1}(u)}{\la_i(u)}&=\frac{P_i(u+1)}{P_i(u)}\qquad\text{for}\quad i=2,\dots,n,
\label{spao}\\
\intertext{and}
\frac{\la_{n+2}(u)}{\la_{n+1}(u)}&=\frac{P_{n+1}(u+2)}{P_{n+1}(u)}.
\label{spala}
\end{align}
The finite-dimensional irreducible
representations of the Yangian $\Y(\osp_{2|2n})$ are in a one-to-one correspondence with
the tuples $\big(\ol Q(u),Q(u),P_2(u),\dots,P_{n+1}(u)\big)$ of monic
polynomials in $u$, where the polynomials
$\ol Q(u)$ and $Q(u)$ are of the same degree and have no common roots.
\end{mthm}

We will refer to the monic
polynomials occurring in these conditions as
the {\em Drinfeld polynomials} of the representation; cf. \cite{d:nr}.

The Appendix is devoted to the polynomial evaluation modules over the Yangian $\Y(\gl_{m|n})$.
Such modules over the Yangian $\Y(\gl_n)$
proved to be useful in the analysis of the representations of the orthosymplectic Yangians;
see Proposition~\ref{prop:xikdp} and also \cite{amr:rp}, \cite{m:ry}.
Many applications of the evaluation homomorphism rely on its $R$-matrix interpretation
going back to Cherednik~\cite{c:ni}. In this interpretation, the polynomial evaluation modules
arise from the Schur--Weyl duality and are constructed as submodules of
the tensor products of the vector representation of the general linear
Lie algebra. As another ingredient, the construction uses
some key steps of the {\em fusion procedure} for the symmetric group
originated in the work of Jucys~\cite{j:yo}. It was
rediscovered in \cite{c:ni}, with detailed arguments given in \cite{n:yc}; see also \cite[Sec.~6.4]{m:yc}.

A super-version of the Schur--Weyl duality involving the general linear Lie superalgebra and
the symmetric group is due to Sergeev~\cite{s:ta} and Berele and Regev~\cite{br:hy}.
It is called the {\em Schur--Sergeev duality} in \cite[Sec.~3.2]{cw:dr}.
We will use this duality to
give an $R$-matrix construction of the polynomial evaluation modules over the Yangian
associated with the Lie superalgebra $\gl_{m|n}$. More precisely,
for any Young diagram $\lambda$ contained in the $(m,n)$-hook we consider
a standard $\la$-tableau $\Uc$ and the associated primitive idempotent $e_{\Uc}\in \CC\Sym_d$,
where $d$ is the number of boxes in $\la$.
By using the action of the symmetric group $\Sym_d$ on the tensor product of the $\ZZ_2$-graded
vector spaces $\CC^{m|n}$, we show that the
subspace $e_{\Uc}(\CC^{m|n})^{\ot d}$ is a representation of
the Yangian $\Y(\gl_{m|n})$. We derive that both this representation
and its twisted version are isomorphic to
evaluation modules over the Yangian and identify them by calculating
the highest weights.

Note that more general
{\em skew representations} of the Yangian $\Y(\gl_{m|n})$ were investigated by Lu and Mukhin~\cite{lm:jt}
with the use of the Schur--Sergeev duality and a super version of the Drinfeld functor.
Those results provide an extension of the corresponding constructions of
Nazarov~\cite{n:rt} for representations of the Yangian $\Y(\gl_{n})$.

\section{Representations of the orthosymplectic Yangian}
\label{sec:ns}

\subsection{Definition and basic properties of the Yangian}
\label{subsec:db}

Fix an integer $n\geqslant 1$ and introduce the
involution $i\mapsto i\pr=2n-i+3$ on
the set $\{1,2,\dots,2n+2\}$.
Consider the $\ZZ_2$-graded vector space $\CC^{2|2n}$ over $\CC$ with the
canonical basis
$e_1,e_2,\dots,e_{2n+2}$, where
the vector $e_i$ has the parity
$\bi\mod 2$ and
\ben
\bi=\begin{cases} 0\qquad\text{for}\quad i=1,1',\\
1\qquad\text{for}\quad i=2,3,\dots,2\pr.
\end{cases}
\een
The endomorphism algebra $\End\CC^{2|2n}$ gets a $\ZZ_2$-gradation with
the parity of the matrix unit $e_{ij}$ found by
$\bi+\bj\mod 2$. We will identify
the algebra of
even matrices over a superalgebra $\Ac$ with the tensor product algebra
$\End\CC^{2|2n}\ot\Ac$, so that a matrix $A=[a_{ij}]$ is regarded as the element
\ben
A=\sum_{i,j=1}^{2n+2}e_{ij}\ot a_{ij}(-1)^{\bi\tss\bj+\bj}\in \End\CC^{2|2n}\ot\Ac.
\een
The involutive matrix {\em super-transposition} $t$ is defined by
$(A^t)_{ij}=A_{j'i'}(-1)^{\bi\bj+\bj}\tss\ta_i\ta_j$,
where we set
\ben
\ta_i=\begin{cases} \phantom{-}1\qquad\text{for}\quad i=1,\dots,n+1,1',\\
-1\qquad\text{for}\quad i=n+2,\dots,2\pr.
\end{cases}
\een
This super-transposition is associated with the bilinear form on the space $\CC^{2|2n}$
defined by the anti-diagonal matrix $G=[\de_{ij'}\tss\ta_i]$.
We will also regard $t$ as the linear map
\beql{suptra}
t:\End\CC^{2|2n}\to \End\CC^{2|2n}, \qquad
e_{ij}\mapsto e_{j'i'}(-1)^{\bi\bj+\bi}\tss\ta_i\ta_j.
\eeq

A standard basis of the general linear Lie superalgebra $\gl_{2|2n}$ is formed by elements $E_{ij}$
of the parity $\bi+\bj\mod 2$ for $1\leqslant i,j\leqslant 2n+2$ with the commutation relations
\ben
[E_{ij},E_{kl}]
=\de_{kj}\ts E_{i\tss l}-\de_{i\tss l}\ts E_{kj}(-1)^{(\bi+\bj)(\bk+\bl)}.
\een
We will regard the orthosymplectic Lie superalgebra $\osp_{2|2n}$
associated with the bilinear from defined by $G$ as the subalgebra
of $\gl_{2|2n}$ spanned by the elements
\ben
F_{ij}=E_{ij}-E_{j'i'}(-1)^{\bi\tss\bj+\bi}\ts\ta_i\ta_j.
\een
The symplectic Lie algebra $\osp_{0|2n}\cong \spa_{2n}$ will be considered as the subalgebra
of $\osp_{2|2n}$ spanned by the elements $F_{ij}$ with $2\leqslant i,j\leqslant 2'$.

Introduce the permutation operator $P$ by
\ben
P=\sum_{i,j=1}^{2n+2} e_{ij}\ot e_{ji}(-1)^{\bj}\in \End\CC^{2|2n}\ot\End\CC^{2|2n}
\een
and set
\ben
Q=\sum_{i,j=1}^{2n+2} e_{ij}\ot e_{i'j'}(-1)^{\bi\bj}\ts\ta_i\ta_j
\in \End\CC^{2|2n}\ot\End\CC^{2|2n}.
\een
The $R$-{\em matrix} associated with $\osp_{2|2n}$ is the
rational function in $u$ given by
\beql{ru}
R(u)=1-\frac{P}{u}+\frac{Q}{u-\ka},\qquad \ka=-n.
\eeq
This is a super-version of the $R$-matrix
originally found in \cite{zz:rf}.
Following \cite{aacfr:rp}, we
define the {\it extended Yangian\/}
$\X(\osp_{2|2n})$
as a $\ZZ_2$-graded algebra with generators
$t_{ij}^{(r)}$ of parity $\bi+\bj\mod 2$, where $1\leqslant i,j\leqslant 2n+2$ and $r=1,2,\dots$,
satisfying certain quadratic relations. To write them down,
introduce the formal series
\beql{tiju}
t_{ij}(u)=\de_{ij}+\sum_{r=1}^{\infty}t_{ij}^{(r)}\ts u^{-r}
\in\X(\osp_{2|2n})[[u^{-1}]]
\eeq
and combine them into the matrix $T(u)=[t_{ij}(u)]$.
Consider the elements of the tensor product algebra
$\End\CC^{2|2n}\ot\End\CC^{2|2n}\ot \X(\osp_{2|2n})[[u^{-1}]]$ given by
\ben
T_1(u)=\sum_{i,j=1}^{2n+2} e_{ij}\ot 1\ot t_{ij}(u)(-1)^{\bi\tss\bj+\bj}\Fand
T_2(u)=\sum_{i,j=1}^{2n+2} 1\ot e_{ij}\ot t_{ij}(u)(-1)^{\bi\tss\bj+\bj}.
\een
The defining relations for the algebra $\X(\osp_{2|2n})$ take
the form of the $RTT$-{\em relation}
\beql{RTT}
R(u-v)\ts T_1(u)\ts T_2(v)=T_2(v)\ts T_1(u)\ts R(u-v).
\eeq
As shown in \cite{aacfr:rp}, the product $T(u-\ka)\ts T^{\tss t}(u)$ is a scalar matrix with
\beql{ttra}
T(u-\ka)\ts T^{\tss t}(u)=c(u)\tss 1,
\eeq
where $c(u)$ is a series in $u^{-1}$. As with the Lie algebra case considered in \cite{amr:rp},
all its coefficients belong to
the center $\ZX(\osp_{2|2n})$ of $\X(\osp_{2|2n})$ and generate the center.

We will also use the {\em extended Yangian} $\X(\osp_{0|2n})$ which is defined
in the same way as $\X(\osp_{2|2n})$, by using the subspace $\CC^{0|2n}\subset\CC^{2|2n}$
with odd basis vectors $e_2,e_3,\dots,e_{2'}$, instead of $\CC^{2|2n}$
and the corresponding $RTT$-relation \eqref{RTT}.
The $R$-matrix is defined
by \eqref{ru} with the modified expressions
\ben
P=-\sum_{i,j=2}^{2n+1} e_{ij}\ot e_{ji}
\Fand
Q=-\sum_{i,j=2}^{2n+1} e_{ij}\ot e_{i'j'}\ts\ta_i\ta_j
\een
and with the value $\ka=-n-1$. It will be convenient to use
the index set $\{2,3,\dots,2'\}$ to label the corresponding generating series
which we will denote by $\bar t_{ij}(u)$. We have an isomorphism
$\X(\osp_{0|2n})\cong \X(\spa_{2n})$ so that the series $\bar t_{ij}(-u)$
satisfy the defining relations for the extended Yangian
$\X(\spa_{2n})$; see \cite{aacfr:rp},
\cite{amr:rp}.

The {\em Yangian} $\Y(\osp_{2|2n})$
is defined as the subalgebra of
$\X(\osp_{2|2n})$ which
consists of the elements stable under
the automorphisms
\beql{muf}
t_{ij}(u)\mapsto f(u)\ts t_{ij}(u)
\eeq
for all series
$f(u)\in 1+u^{-1}\CC[[u^{-1}]]$.
We have the tensor product decomposition
\beql{tensordecom}
\X(\osp_{2|2n})=\ZX(\osp_{2|2n})\ot \Y(\osp_{2|2n}).
\eeq
The Yangian $\Y(\osp_{2|2n})$ is isomorphic to the quotient
of $\X(\osp_{2|2n})$
by the relation $c(u)=1$.

The defining relations \eqref{RTT} can be written
with the use of super-commutator in terms of the series \eqref{tiju} as follows:
\begin{align}
\big[\tss t_{ij}(u),t_{kl}(v)\big]&=\frac{1}{u-v}
\big(t_{kj}(u)\ts t_{il}(v)-t_{kj}(v)\ts t_{il}(u)\big)
(-1)^{\bi\tss\bj+\bi\tss\bk+\bj\tss\bk}
\non\\
{}&-\frac{1}{u-v-\ka}
\Big(\de_{k i\pr}\sum_{p=1}^{2n+2}\ts t_{pj}(u)\ts t_{p'l}(v)
(-1)^{\bi+\bi\tss\bj+\bj\tss\bp}\ts\ta_i\ta_p
\label{defrel}\\
&\qquad\qquad\qquad
{}-\de_{l j\pr}\sum_{p=1}^{2n+2}\ts t_{k\tss p'}(v)\ts t_{ip}(u)
(-1)^{\bj+\bp+\bi\tss\bk+\bj\tss\bk+\bi\tss\bp}\ts\ta_j\ta_p\Big).
\non
\end{align}
Note that the mapping
\beql{shift}
t_{ij}(u)\mapsto t_{ij}(u+a),\quad a\in \CC,
\eeq
defines an automorphism of $\X(\osp_{2|2n})$.

The universal enveloping algebra $\U(\osp_{2|2n})$ can be regarded as a subalgebra of
$\X(\osp_{2|2n})$ via the embedding
\beql{emb}
F_{ij}\mapsto \frac12\big(t_{ij}^{(1)}-t_{j'i'}^{(1)}(-1)^{\bj+\bi\bj}\ts\ta_i\ta_j\big)(-1)^{\bi}.
\eeq
This fact relies on the Poincar\'e--Birkhoff--Witt theorem for the orthosymplectic Yangian
which was pointed out in \cite{aacfr:rp}.
It states that the associated graded algebra
for $\Y(\osp_{2|2n})$ is isomorphic to $\U(\osp_{2|2n}[u])$.
This implies that the algebra $\X(\osp_{2|2n})$ is generated by
the coefficients of the series $c(u)$ and $t_{ij}(u)$ with $i+j\leqslant 2n+3$,
excluding $t_{11'}(u)$ and $t_{1'1}(u)$. Moreover, given any total ordering
on the set of the generators, the ordered monomials with the powers of odd generators
not exceeding $1$, form a basis of the algebra.
The theorem follows by using essentially the same
arguments as in the non-super case; see \cite[Sec.~3]{amr:rp}.

The extended Yangian $\X(\osp_{2|2n})$ is a Hopf algebra with the coproduct
defined by
\beql{Delta}
\De: t_{ij}(u)\mapsto \sum_{k=1}^{2n+2} t_{ik}(u)\ot t_{kj}(u).
\eeq
The coproduct on the algebra $\X(\osp_{0|2n})$ is defined by the same formula
for the series $\bar t_{ij}(u)$ instead of $t_{ij}(u)$, with the sum taken
over $k=2,3,\dots,2n+1$.

\subsection{Highest weight representations}
\label{subsec:hw}

The proof of the necessity of the conditions of the Main Theorem will rely on the following
reduction property for representations of the extended Yangians $\X(\osp_{2|2n})$
which is verified in the same way as for the case of $\X(\osp_{1|2n})$ in \cite[Prop.~4.1]{m:ry};
cf. \cite[Lemma~5.13]{amr:rp}.
For an arbitrary $\X(\osp_{2|2n})$-module $V$ set
\beql{vplus}
V^+=\{\eta\in V\ |\ t_{1j}(u)\ts\eta=0\quad\text{for}\quad j>1\fand
t_{i\tss 1'}(u)\ts\eta=0\quad\text{for}\quad i<1'\}.
\eeq

\bpr\label{prop:vplus}
The subspace $V^+$ is stable under the action of the operators $t_{ij}(u)$
subject to $2\leqslant i,j\leqslant 2'$.
Moreover, the assignment $\bar t_{ij}(u)\mapsto t_{ij}(u)$
defines
a representation of the extended Yangian $\X(\osp_{0|2n})$ on $V^+$.
\qed
\epr

A representation $V$ of the algebra $\X(\osp_{2|2n})$
is called a {\em highest weight representation}
if there exists a nonzero vector
$\xi\in V$ such that $V$ is generated by $\xi$,
\begin{alignat}{2}
t_{ij}(u)\ts\xi&=0 \qquad &&\text{for}
\quad 1\leqslant i<j\leqslant 1', \qquad \text{and}\non\\
t_{ii}(u)\ts\xi&=\la_i(u)\ts\xi \qquad &&\text{for}
\quad i=1,\dots,1',
\label{trianb}
\end{alignat}
for some formal series
\beql{laiu}
\la_i(u)\in 1+u^{-1}\CC[[u^{-1}]].
\eeq
The vector $\xi$ is called the {\em highest vector}
of $V$.

\bpr\label{prop:nontrvm}
The series $\la_i(u)$ associated with a highest weight representation $V$
satisfy
the consistency conditions
\beql{nontrvm}
\la_i(u)\tss \la_{i\pr}(u+n-i+2)=\la_{i+1}(u)\tss \la_{(i+1)'}(u+n-i+2)
\eeq
for $i=1,\dots,n$.
Moreover, the coefficients of the series $c(u)$ act in the representation
$V$ as the multiplications by scalars
determined by
$
c(u)\mapsto \la_1(u)\tss \la_{1'}(u+n).
$
\epr

\bpf
Introduce the subspace $V^+$ by \eqref{vplus} and note that
the vector $\xi$ belongs to $V^+$. By applying Proposition~\ref{prop:vplus} we find that
the cyclic span
$\X(\osp_{0|2n})\ts\xi$ is a highest weight submodule over $\X(\osp_{0|2n})$ associated with the series
$\la_2(u),\dots,\la_{2'}(u)$.
Due to the isomorphism
$\X(\osp_{0|2n})\cong \X(\spa_{2n})$, conditions
\eqref{nontrvm} for $i=2,\dots,n$ follow from the consistency conditions
for the highest weight modules over $\X(\spa_{2n})$; see \cite[Prop.~5.14]{amr:rp}.
Furthermore, using the defining relations \eqref{defrel}, we get
\ben
t_{12}(u)\ts t_{1'2'}(v)\ts \xi=
-\frac{1}{u-v-\kappa}\ts\Big(t_{12}(u)\ts t_{1'2'}(v)+
\la_{1}(u)\ts \la_{1'}(v)-\la_{2}(u)\ts \la_{2'}(v)\Big)\ts \xi
\een
and so
\ben
(u-v-\kappa+1)\ts t_{12}(u)\ts t_{1'2'}(v)\ts \xi
=\big({-}\la_{1}(u)\ts \la_{1'}(v)+\la_{2}(u)\ts \la_{2'}(v)\big)\ts \xi.
\een
Setting $v=u-\kappa+1=u+n+1$, we obtain
\eqref{nontrvm} for $i=1$. Finally, the last part of the proposition
is obtained by using the expression for $c(u)$ implied by taking the $(1',1')$ entry in
the matrix relation \eqref{ttra}.
\epf

By Proposition~\ref{prop:nontrvm}, the series $\la_i(u)$ in \eqref{trianb}
with $i>n+2$ are uniquely
determined by the first $n+2$ series. The corresponding $(n+2)$-tuple
$\la(u)=(\la_{1}(u),\dots,\la_{n+2}(u))$ will be called
the {\em highest weight} of $V$.

For an arbitrary $(n+2)$-tuple $\la(u)=(\la_{1}(u),\dots,\la_{n+2}(u))$
of formal series of the form \eqref{laiu}
define
the {\em Verma module} $M(\la(u))$ as the quotient of the algebra $\X(\osp_{2|2n})$ by
the left ideal generated by all coefficients of the series $t_{ij}(u)$
with $1\leqslant i<j\leqslant 2n+2$, and $t_{ii}(u)-\la_i(u)$ for
$i=1,\dots,2n+2$, assuming that the series $\la_i(u)$
with $i=n+3,\dots,2n+2$ are defined to satisfy the consistency conditions \eqref{nontrvm}.

The Poincar\'e--Birkhoff--Witt theorem for the algebra $\X(\osp_{2|2n})$
implies that the Verma module $M(\la(u))$
is nonzero, and we denote by $L(\la(u))$ its irreducible quotient.
It is clear that the isomorphism
class of $L(\la(u))$ is determined by $\la(u)$.

The first part of the Main Theorem is implied by the following proposition
whose proof is the same as in the non-super case; see \cite[Thm~5.1]{amr:rp}.

\bpr\label{prop:fdhw}
Every finite-dimensional irreducible representation of the algebra $\X(\osp_{2|2n})$
is a highest weight representation.
Moreover, it contains a unique, up to a constant factor, highest vector.
\qed
\epr

We will now suppose that the representation $L(\la(u))$ of $\X(\osp_{2|2n})$
is finite-dimensional and derive the conditions on the highest weight $\la(u)$ given
in the Main Theorem. First consider
the subalgebra $\Y_0$ of $\X(\osp_{2|2n})$ generated by the coefficients
of the series $t_{11}(u)$, $t_{12}(u)$, $t_{21}(u)$ and $t_{22}(u)$.
This subalgebra is isomorphic to the Yangian $\Y(\gl_{1|1})$, and the cyclic span
$\Y_0\tss \xi$ of the highest vector $\xi$ is a finite-dimensional
module over $\Y_0$ with the highest weight $(\la_1(u),\la_2(u))$.
This implies that condition \eqref{gloo} must hold by \cite[Thm.~4]{zh:rs};
see also \cite[Prop.~5.1]{m:ry}.

Furthermore, Proposition~\ref{prop:vplus} implies that the subspace
$L(\la(u))^+$ is a module over the extended Yangian
$\X(\osp_{0|2n})\cong \X(\spa_{2n})$. The vector $\xi$ generates a highest weight
$\X(\spa_{2n})$-module
with the highest weight $(\la_{2}(-u),\dots,\la_{n+2}(-u))$.
Since this module is finite-dimensional, conditions
\eqref{spao} and \eqref{spala} now follow by
\cite[Thm.~5.16]{amr:rp}.

\subsection{Fundamental representations}
\label{sec:fm}

Our next step in the proof of the Main Theorem is to show that the conditions
\eqref{gloo}, \eqref{spao} and \eqref{spala} are sufficient for the representation
$L(\la(u))$ of $\X(\osp_{2|2n})$ to be finite-dimensional.
The tuple
\beql{tdp}
\big(\ol Q(u),Q(u),P_2(u),\dots,P_{n+1}(u)\big)
\eeq
of Drinfeld
polynomials determines the highest weight
$\la(u)$ of the representation up to a simultaneous multiplication of all components
$\la_i(u)$ by a formal series $f(u)\in 1+u^{-1}\CC[[u^{-1}]]$.
This operation corresponds to twisting the action of
the algebra $\X(\osp_{2|2n})$ on $L(\la(u))$ by the automorphism \eqref{muf}.
Hence, it suffices to prove that
a particular module $L(\la(u))$ corresponding to a given tuple \eqref{tdp}
is finite-dimensional.

Suppose that $L(\nu(u))$ and $L(\mu(u))$ are the irreducible
highest weight modules with the highest weights
\beql{munu}
\nu(u)=\big(\nu_{1}(u),\dots,\nu_{n+2}(u)\big)\Fand
\mu(u)=\big(\mu_{1}(u),\dots,\mu_{n+2}(u)\big).
\eeq
By the coproduct rule
\eqref{Delta},
the cyclic span $\X(\osp_{2|2n})(\xi\ot\xi')$ of the tensor product
of the respective highest vectors $\xi\in L(\nu(u))$ and $\xi'\in L(\mu(u))$
is a highest weight
module
with the highest weight
\beql{munup}
\big(\nu_{1}(u)\tss\mu_{1}(u),\dots,\nu_{n+2}(u)\tss\mu_{n+2}(u)\big).
\eeq
This observation implies the corresponding transition rule for the associated tuples
of Drinfeld polynomials. Namely, if
the highest weights in \eqref{munu}
are associated with
the respective tuples
\ben
\big(\ol Q(u),Q(u),P_2(u),\dots,P_{n+1}(u)\big)\Fand
\big(\ol Q^{\circ}(u),Q^{\circ}(u),P^{\circ}_2(u),\dots,P^{\circ}_{n+1}(u)\big),
\een
then
the highest weight \eqref{munup} is associated with the tuple
\beql{traru}
\Big(\frac{\ol Q(u)\tss\ol Q^{\circ}(u)}{d(u)},\frac{Q(u)\tss Q^{\circ}(u)}{d(u)},
P^{}_2(u)P^{\circ}_2(u),\dots,P^{}_{n+1}(u)P^{\circ}_{n+1}(u)\Big),
\eeq
where $d(u)$ is the monic polynomial in $u$ defined as the greatest common divisor of
the polynomials $\ol Q(u)\tss\ol Q^{\circ}(u)$ and $Q(u)\tss Q^{\circ}(u)$.
Therefore, it will be enough
to show that the {\em fundamental representations} of $\X(\osp_{2|2n})$
are finite-dimensional. They correspond to the tuples of Drinfeld polynomials
of the form
\beql{fundam}
\big(u+\al,u+\be,1,\dots,1\big)\Fand
\big(1,\dots,1,u+\ga,1,\dots,1\big),
\eeq
for $\al\ne\be$, where $u+\ga$ represents the polynomial $P_d(u)$ and $d$ runs
over the values $2,\dots,n+1$.
Note that by applying the shift automorphism \eqref{shift},
we can assume that $\be=0$ in the first tuple and use a particular value of $\ga$
in the second tuple.

We will begin with the first tuple in \eqref{fundam}.
The required property is implied by the following.

\bpr\label{prop:laon}
The representation $L(\la(u))$ of the extended Yangian $\X(\osp_{2|2n})$
with the components of the
highest weight
given by $\la_1(u)=1+\al\tss u^{-1}$
and $\la_i(u)=1$ for $i=2,3,\dots,n+2$, is finite-dimensional.
\epr

\bpf
Let $\xi$ denote the highest vector of $L(\la(u))$.
We will show first that all vectors $t_{kl}^{(r)}\xi$ with $k,l\in\{2,3,\dots,2'\}$
and $r\geqslant 1$ are zero. This holds by definition for $k\leqslant l$, while
the defining relations \eqref{defrel} imply that for $k>l$
we have the properties $t_{1j}(u)\tss t_{kl}^{(r)}\xi=0$ for $j>1$ and
$t_{i1'}(u)\tss t_{kl}^{(r)}\xi=0$ for $i<1'$. This means that the vector
$t_{kl}^{(r)}\xi$ belongs to the subspace $L(\la(u))^+$ defined in \eqref{vplus}.
By Proposition~\ref{prop:vplus}, the subspace $L(\la(u))^+$
is a representation of the algebra $\X(\osp_{0|2n})$, and $\xi$ generates
a highest weight module over this algebra with the highest weight $(1,\dots,1)$.
The irreducible quotient of this module is one-dimensional, which implies that
the vector $t_{kl}^{(r)}\xi\in L(\la(u))$ is annihilated
by all coefficients of the series $t_{ij}(u)$ with $1\leqslant i<j\leqslant 1'$
and therefore must be zero.

We now prove two lemmas providing
the formulas for the action of some other elements
of the extended Yangian on the highest vector $\xi$.

\ble\label{lem:vani}
We have $t_{k1}^{(r)}\xi=0$
for all $k=2,3,\dots,2'$ and $r\geqslant 2$.
\ele

\bpf
As in the above argument, we will check
that the vectors $t_{k1}^{(r)}\xi$ are annihilated by all
coefficients of the series $t_{ij}(u)$ with $1\leqslant i<j\leqslant 1'$.
The Poincar\'e--Birkhoff--Witt theorem and the defining relations
imply that it is sufficient to check this property
for the coefficients of the series $t_{i,i+1}(u)$ with $i=1,\dots,n+1$.
Suppose first that $2\leqslant k\leqslant n+1$. Then \eqref{defrel} gives
\ben
\big[\tss t_{k1}(v),t_{i,i+1}(u)\big]=\frac{1}{v-u}
\big(t_{i1}(v)\ts t_{k,i+1}(u)-t_{i1}(u)\ts t_{k,i+1}(v)\big)(-1)^{\bi}.
\een
Since $t_{k,i+1}(u)\tss\xi=\de_{k,i+1}\tss\xi$, we derive
\beql{tiipo}
t_{i,i+1}(u)\tss t_{k1}(v)\tss\xi=\frac{1}{v-u}\ts \de_{k,i+1}\tss
\big(t_{i1}(v)-t_{i1}(u)\big)\tss\xi.
\eeq
Now proceed by induction on $k$.
For $k=2$, by \eqref{tiipo} we have
\ben
t_{12}(u)\tss t_{21}(v)\tss\xi=\frac{1}{v-u}\ts \big(t_{11}(v)-t_{11}(u)\big)\tss\xi
=-\al\tss u^{-1}\tss v^{-1}\tss\xi
\een
which implies
that $t_{21}(v)\tss\xi=v^{-1}\tss t_{21}^{(1)}\tss\xi$. Using \eqref{tiipo}
for the induction step, we can conclude that $t_{k1}(v)\tss\xi=v^{-1}\tss t_{k1}^{(1)}\tss\xi$
for all $2\leqslant k\leqslant n+1$.

Now let $(n+1)'\leqslant k\leqslant 2'$ and use induction on $k$ starting from $k=(n+1)'$.
Note that \eqref{tiipo} holds for these values of $k$ unless $i=k'$. By applying \eqref{defrel}
to the super-commutator $\big[t_{(n+1)',1}(v),t_{n+1,(n+1)'}(u)\big]$ and using the assumptions
$\la_i(u)=1$ for $i=2,3,\dots,n+2$,
we get
\ben
\bal
t_{n+1,(n+1)'}(u)\ts t_{(n+1)',1}(v)\ts \xi
&=\frac{1}{v-u}\ts \big(t_{n+1,1}(v)-t_{n+1,1}(u)\big)\tss\xi\\[0.4em]
{}&+\frac{1}{v-u+n}\ts t_{11}(v)\ts t_{1',(n+1)'}(u)\tss\xi
+\frac{1}{v-u+n}\ts t_{n+1,1}(v)\tss\xi.
\eal
\een
By another application of \eqref{defrel} we derive
\ben
\frac{v-u+n+1}{v-u+n}\ts t_{11}(v)\ts t_{1',(n+1)'}(u)\tss\xi=
\frac{v+\al}{v}\ts t_{1',(n+1)'}(u)\tss\xi-\frac{1}{v-u+n}\ts t_{n+1,1}(v)\tss\xi.
\een
Therefore, the previous expression can be written in the form
\begin{align}
t_{n+1,(n+1)'}(u)\ts t_{(n+1)',1}(v)\ts \xi
&=\frac{1}{v-u}\ts \big(t_{n+1,1}(v)-t_{n+1,1}(u)\big)\tss\xi
\label{tnpn}\\[0.4em]
{}&+\frac{v+\al}{v(v-u+n+1)}\ts t_{1',(n+1)'}(u)\tss\xi
+\frac{1}{v-u+n+1}\ts t_{n+1,1}(v)\tss\xi.
\non
\end{align}
On the other hand, by calculating the super-commutator $[t_{n+1,1}(v),t_{1'1'}(u)]$
we obtain
\ben
t_{1'1'}(u)\ts t_{n+1,1}(v)\ts\xi=\frac{v-u+n}{v-u+n+1}\ts t_{n+1,1}(v)\ts t_{1'1'}(u)\ts\xi
-\frac{1}{v-u+n+1}\ts t_{1',(n+1)'}(u)\tss\xi.
\een
By the consistency conditions of Proposition~\ref{prop:nontrvm}, we find that
\ben
t_{1'1'}(u)\ts\xi=\frac{u-n-1}{u+\al-n-1}\ts\xi,
\een
and so the residue at $v=u-n-1$ gives
\ben
t_{1',(n+1)'}(u)\tss\xi=-\frac{u-n-1}{u+\al-n-1}\ts t_{n+1,1}(u-n-1)\ts\xi.
\een
As we showed above, $t_{n+1,1}(v)\tss\xi=v^{-1}\tss t_{n+1,1}^{(1)}\tss\xi$,
so that \eqref{tnpn} simplifies to
\ben
t_{n+1,(n+1)'}(u)\ts t_{(n+1)',1}(v)\ts \xi=-\Big(\frac{1}{u\tss v}+\frac{1}{(u+\al-n-1)\tss v}\Big)
\ts t_{n+1,1}^{(1)}\tss\xi.
\een
Hence we may conclude that the vectors
$t_{(n+1)',1}^{(r)}\xi$ with $r\geqslant 2$ are annihilated by all
coefficients of the series $t_{i,i+1}(u)$ with $i=1,\dots,n+1$
and so $t_{(n+1)',1}^{(r)}\xi=0$.

Now suppose that $k=l'$ with $2\leqslant l\leqslant n$. Calculating as above, we find that
\beql{tlpo}
t_{l,l+1}(u)\tss t_{l',1}(v)\ts \xi=\frac{1}{v-u+n}\ts\big(t_{11}(v)\ts t_{1',l+1}(u)\ts\xi
-t_{(l+1)',1}(v)\ts\xi\big).
\eeq
Furthermore, using the defining relations again, we derive the relation
\ben
t_{11}(v)\ts t_{1',l+1}(u)\ts\xi=\frac{1}{v-u+n+1}\ts t_{(l+1)',1}(v)\ts\xi
+\frac{(v+\al)(v-u+n)}{v\tss(v-u+n+1)}\ts t_{1',l+1}(u)\ts\xi,
\een
which by taking the residue at $v=u-n-1$ also implies
\ben
t_{(l+1)',1}(u-n-1)\ts\xi=\frac{u+\al-n-1}{u-n-1}\ts t_{1',l+1}(u)\ts\xi.
\een
Since
$
t_{(l+1)',1}(v)\ts\xi=v^{-1}\ts t_{(l+1)',1}^{(1)}\ts\xi
$
by the induction hypothesis,
relation \eqref{tlpo} simplifies to
\ben
t_{l,l+1}(u)\tss t_{l',1}(v)\ts \xi=
\frac{1}{(u+\al-n-1)\tss v}
\ts t_{(l+1)',1}^{(1)}\ts\xi,
\een
thus proving that the vectors
$t_{l',1}^{(r)}\xi$ with $r\geqslant 2$ are annihilated by the
coefficients of the series $t_{l,l+1}(u)$.
The proof is competed by taking into account \eqref{tiipo}
for the remaining values $i\ne l$.
\epf

\ble\label{lem:for}
We have
\ben
t_{1'1}(v)\ts\xi=\frac{1}{v\ts(v+\al-n-1)}\ts \sum_{k=2}^{n+1}\tss t_{k'1}^{(1)}\tss t_{k\tss 1}^{(1)}\ts\xi.
\een
\ele

\bpf
Calculate the commutator $[t_{11}(v),t_{1'1}(u)]$ by \eqref{defrel}
and rearrange the terms to get
\ben
\bal
\frac{v-u+n+1}{v-u+n}\ts t_{11}(v)\ts t_{1'1}(u)&=
\frac{n}{(v-u)(v-u+n)}\ts t_{1'1}(v)\ts t_{11}(u)\\
{}&+\frac{v-u-1}{v-u}\ts t_{1'1}(u)\ts t_{11}(v)
-\frac{1}{v-u+n}\ts\sum_{k=2}^{2'}t_{k1}(v)\ts t_{k'1}(u)\ts\ta_k.
\eal
\een
By setting $u=v+n+1$ we come to the relation
\ben
\frac{n}{n+1}\ts t_{1'1}(v)\ts t_{11}(v+n+1)
+\frac{n+2}{n+1}\ts t_{1'1}(v+n+1)\ts t_{11}(v)
=-\sum_{k=2}^{2'}t_{k1}(v)\ts t_{k'1}(v+n+1)\ts\ta_k.
\een
Apply both sides of the relation to the highest vector $\xi$.
By Lemma~\ref{prop:laon} we have
\ben
t_{k'1}(v+n+1)\ts\xi=(v+n+1)^{-1}\ts t_{k'1}^{(1)}\ts\xi.
\een
Hence, using the super-commutator
$
[t_{k1}(v),t_{k'1}^{(1)}]=-t_{1'1}(v)\ts\ta_k
$
we obtain
\begin{multline}
\frac{n(v+\al-n-1)}{(n+1)(v+n+1)}\ts t_{1'1}(v)\ts \xi
+\frac{(n+2)(v+\al)}{(n+1)v}\ts t_{1'1}(v+n+1)\ts \xi\\
{}=\frac{1}{v\ts(v+n+1)}\ts \sum_{k=2}^{2'}\tss \ta_k\ts t_{k'1}^{(1)}
\tss t_{k\tss 1}^{(1)}\ts\xi.
\non
\end{multline}
This relation uniquely determines the series $t_{1'1}(v)\ts \xi$. It is easily seen to be
satisfied by
\ben
t_{1'1}(v)\ts \xi=\frac{1}{2\tss v\ts(v+\al-n-1)}\ts
\sum_{k=2}^{2'}\tss \ta_k\ts t_{k'1}^{(1)}\tss t_{k\tss 1}^{(1)}\ts\xi.
\een
In particular, $t_{1'1}^{(1)}\ts \xi=0$ and so $[t_{k'1}^{(1)}, t_{k\tss 1}^{(1)}]\ts\xi=0$
yielding the required relation.
\epf

Returning to the proof of the proposition, we can now derive that
the representation
$L(\la(u))$ is spanned by the vectors of the form
\beql{basve}
t_{k_1 1}^{(1)}\dots t_{k_s 1}^{(1)}\ts\xi,\qquad 2\leqslant k_1<\dots<k_s\leqslant 2',
\eeq
with $s\geqslant 0$.
Indeed, by the first part of the proof,
and the Poincar\'e--Birkhoff--Witt theorem applied to a suitable ordering of the generators,
it is enough to verify that
the span of these vectors is stable under the action of all generators of the form
$t_{k\tss 1}^{(r)}$ with $k=2,\dots,2'$. Note that by taking the $(1',1)$ entry
in the matrix relation \eqref{ttra},
we find that $t_{1'1}^{(1)}=0$ in the algebra $\X(\osp_{2|2n})$. This implies that
the elements $t_{k\tss 1}^{(1)}$ with $k=2,\dots,2'$ pairwise anticommute.
They also have the property $[t^{}_{1'1}(v),t_{k\tss 1}^{(1)}]=0$.
Hence, applying the series $t_{k\tss 1}(v)$ to a basis vector \eqref{basve}, we get
\ben
t^{}_{k\tss 1}(v)\ts t_{k_1 1}^{(1)}\dots t_{k_s 1}^{(1)}\ts\xi
=\sum_{i=1}^s(-1)^{i-1}\ts t_{k_1 1}^{(1)}\dots [t^{}_{k 1}(v),t_{k_i\tss 1}^{(1)}]
\dots t_{k_s 1}^{(1)}\ts\xi+(-1)^s\ts t_{k_1 1}^{(1)}\dots t_{k_s 1}^{(1)}\ts
t^{}_{k\tss 1}(v)\ts\xi.
\een
Since $[t^{}_{k\tss 1}(v),t_{k_i 1}^{(1)}]=-\de_{k'k_i}t_{1'1}(v)\tss\ta_k$,
the sum on the right hand side equals
\ben
\ta_k\ts\sum_{i=1}^s(-1)^{i}\ts\de_{k'k_i}\ts t_{k_1 1}^{(1)}\dots \wh t_{k_i 1}^{\ts(1)}
\dots t_{k_s 1}^{(1)}\ts t_{1'1}(v)\ts\xi,
\een
with the hat indicating the factor to be skipped.
By Lemma~\ref{lem:vani} we have $t^{}_{k\tss 1}(v)\ts\xi=v^{-1}\ts t_{k\tss 1}^{(1)}\ts\xi$,
and together with the relation of Lemma~\ref{lem:for} this proves that
the span of the vectors \eqref{basve} is stable under the action of all
operators $t_{k1}^{(r)}$ with $k=2,\dots,2'$. Since the representation
$L(\la(u))$ is spanned by these vectors, its dimension does not exceed $4^n$.
\epf

Let $V(\mu)$ denote the irreducible highest weight representation of the Lie superalgebra
$\osp_{2|2n}$ with the highest weight $\mu=(\mu_1,\dots,\mu_{n+1})$ with respect to the
upper-triangular Borel subalgebra. This means that $F_{ii}\tss\ze=\mu_i\ts\ze$ for $i=1,\dots,n+1$
and $F_{ij}\ts\ze =0$ for $i<j$, where $\ze$ is the highest vector of $V(\mu)$.
The modules of the form
$V(\al,0,\dots,0)$ are known to be {\em typical}, if and only if $\al$ does not
belong to the set $\{0,1,\dots,n-1\}\cup\{n+1,\dots,2n\}$; see \cite{k:rc}.

\bco\label{cor:exls}
The typical representations $V(\al,0,\dots,0)$ of the Lie superalgebra $\osp_{2|2n}$ extend
to modules over the Yangian $\X(\osp_{2|2n})$.
\eco

\bpf
Equip the representation $L(1+\al\ts u^{-1},1,\dots,1)$ of $\X(\osp_{2|2n})$ with the action
of $\osp_{2|2n}$ via the embedding \eqref{emb}. By taking the $(i,j)$
entry in the matrix relation \eqref{ttra} for $i\ne j$ we find
\ben
t_{ij}^{(1)}+t_{j'i'}^{(1)}(-1)^{\bj+\bi\bj}\ts\ta_i\ta_j=0.
\een
Therefore, under the embedding we have $F_{ij}\mapsto t_{ij}^{(1)}(-1)^{\bi}$.

On the other hand, any typical representation $V(\al,0,\dots,0)$
is isomorphic to the corresponding Kac module
and so is equipped with a basis of the form
\ben
F_{k_1 1}\dots F_{k_s 1}\ts\ze,\qquad 2\leqslant k_1<\dots<k_s\leqslant 2',
\een
with $s\geqslant 0$. The proof of Proposition~\ref{prop:laon} shows that the mapping
$\ze\to\xi$ extends to an $\osp_{2|2n}$-module isomorphism
$V(\al,0,\dots,0)\to L(1+\al\ts u^{-1},1,\dots,1)$.
\epf

We will now turn to the second family of
fundamental representations in \eqref{fundam}.
We will show that they can be constructed as subquotients of tensor products
of the vector representations of $\X(\osp_{2|2n})$ and representations
of the first family in \eqref{fundam}.

The {\em vector representation} of $\X(\osp_{2|2n})$ on $\CC^{2|2n}$ is defined by
\beql{vectre}
t_{ij}(u)\mapsto \de_{ij}+u^{-1}\tss
e_{ij}(-1)^{\bi}-(u-n)^{-1}\tss e_{j'i'}(-1)^{\bi\bj}\ts\ta_i\ta_j.
\eeq
The homomorphism property
follows from the $RTT$-relation \eqref{RTT} and
the Yang--Baxter equation satisfied by $R(u)$ as in \eqref{ybe}; cf. \cite{aacfr:rp}.
The mapping $T(u)\mapsto R(u)$ defines an algebra homomorphism
$\X(\osp_{2|2n})\to \End\CC^{2|2n}$, and to get \eqref{vectre} we take
its composition with the
automorphism $t_{ij}(u)\mapsto t_{i'j'}(u+n)\ts\ta_i\ta_j(-1)^{\bi+\bj}$.

Now use the coproduct \eqref{Delta} and the shift automorphism \eqref{shift}
to equip
the tensor square $V=(\CC^{2|2n})^{\ot 2}$
with the action of $\X(\osp_{2|2n})$ by setting
\beql{tijtprei}
t_{ij}(u)\mapsto
\sum_{l=1}^{2n+2} t_{i\tss l}(u-1)\ot t_{lj}(u),
\eeq
where the generators act in the respective copies of the vector space
$\CC^{2|2n}$ via the rule \eqref{vectre}.
Introduce the vectors $w_k\in V$ by $w_k=e_1\ot e_k-e_k\ot e_1$
for $k=2,3,\dots,2'$ and denote by $W$ the subspace of $V$
spanned by these vectors.

The vector representation $\CC^{0|2n}$ of the extended Yangian
$\X(\osp_{0|2n})$ is defined by
the rule similar to \eqref{vectre}; see \cite{amr:rp}:
\beql{vectreo}
\bar t_{ij}(u)\mapsto \de_{ij}-u^{-1}\tss
e_{ij}+(u-n-1)^{-1}\tss e_{j'i'}\ts\ta_i\ta_j,\qquad 2\leqslant i,j\leqslant 2'.
\eeq

\bpr\label{prop:vectv}\quad {\rm (i)}\ \
The subspace $W$ is invariant with respect to
the action of the operators $t_{ij}(u)$
with $2\leqslant i,j\leqslant 2'$.
Moreover, the assignment $\bar t_{ij}(u)\mapsto t_{ij}(u)$
defines
a representation of the algebra $\X(\osp_{0|2n})$ on $W$
isomorphic to the vector representation.
\par
{\rm (ii)}\ \  The subspace $W$ is annihilated by the operators $t_{1j}(u)$ with $j>1$ and
$t_{i\tss 1'}(u)$ with $i<1'$. Moreover, the coefficients of the series $t_{11}(u)$
act on $W$ as multiplications by scalars determined by $t_{11}(u)\mapsto 1+u^{-1}$.
\epr

\bpf
Assuming that $2\leqslant i,j\leqslant 2'$, for
the action of the operator $t_{ij}(u)$
we have
\ben
\bal
t_{ij}(u)\ts w_k&=t_{ij}(u)\ts(e_1\ot e_k-e_k\ot e_1)\\
&{}=\sum_{l=1}^{2n+2} t_{i\tss l}(u-1)\ts e_1\ot t_{lj}(u)\ts e_k
-\sum_{l=1}^{2n+2} (-1)^{l+j}\ts t_{i\tss l}(u-1)\ts e_k\ot t_{lj}(u)\ts e_1.
\eal
\een
Due to \eqref{vectre},
a nonzero contribution can only come from the terms
\begin{multline}
t_{i1}(u-1)\ts e_1\ot t_{1j}(u)\ts e_k+t_{ii}(u-1)\ts e_1\ot t_{ij}(u)\ts e_k\\[0.3em]
-t_{ij}(u-1)\ts e_k\ot t_{jj}(u)\ts e_1+t_{i1'}(u-1)\ts e_k\ot t_{1'j}(u)\ts e_1
\non
\end{multline}
and the calculation yields the formula
\ben
t_{ij}(u)\ts w_k=\de_{ij}\ts w_k-\de_{kj}\ts u^{-1}\ts w_i+\de_{ki'}\ts\ta_i\tss\ta_j\ts
 (u-n-1)^{-1}\ts w_{j'}.
\een
Taking onto account \eqref{vectreo}, we thus prove
the first part of the proposition. The second part is verified by
a similar calculation.
\epf

For any $a\in\CC$ we will denote by $V_a$ the representation
of the algebra $\X(\osp_{2|2n})$ on the space $V=(\CC^{2|2n})^{\ot 2}$
obtained by twisting the action \eqref{tijtprei}
with the automorphism \eqref{shift}. Similarly, we let $W_a$ denote the
representation of $\X(\osp_{0|2n})$ on $W$ twisted by the same shift automorphism
on $\X(\osp_{0|2n})$.
For $d=1,\dots,n$ introduce the vectors
\ben
\xi_d=\sum_{\si\in\Sym'_d} \sgn\si\cdot
w_{\si(2)}\ot w_{\si(3)}\ot\dots\ot w_{\si(d+1)} \in V_0\ot V_{-1}\ot\dots\ot V_{-d+1},
\een
where $\Sym'_d$ denotes the group of permutations of the set $\{2,\dots,d+1\}$.

\bpr\label{prop:xikdp}
The cyclic span $\X(\osp_{2|2n})\ts\xi_d$ is a highest weight module
whose tuple of Drinfeld polynomials \eqref{tdp} has the form
\beql{tudp}
(u+1,u-d,1,\dots,1,u-d,1,\dots,1)\quad \text{with}\quad P_{d+1}(u)=u-d,
\eeq
for $d=1,\dots,n-1$, while for $d=n$
it has the form
\beql{tudpn}
(u+1,u-n,1,\dots,1,u-n-1).
\eeq
\epr

\bpf
It follows from Proposition~\ref{prop:vectv} and the coproduct rule \eqref{Delta}
that for all $2\leqslant i,j\leqslant 2'$ the image $t_{ij}(u)\ts\xi_d$
coincides with $\bar t_{ij}(u)\ts\xi_d$ in the representation
$W_0\ot W_{-1}\ot\dots\ot W_{-d+1}$
of the algebra $\X(\osp_{0|2n})$. Since the series $\bar t_{ij}(-u)$ satisfy
the defining relations of the extended Yangian $\X(\spa_{2n})$, we derive from
the proof of \cite[Thm.~5.16]{amr:rp} (see also Corollary~\ref{cor:fungl} below), that
the vector $\xi_d$ has the properties
\ben
t_{ij}(u)\ts \xi_d=0\qquad\text{for}\quad 2\leqslant i<j\leqslant 2'
\een
and
\beql{tiixi}
t_{ii}(u)\ts \xi_d=\begin{cases}
\dfrac{u-d}{u-d+1}\ts\xi_d \qquad&\text{for}
\quad i=2,\dots, d+1,\\[0.6em]
\ \xi_d\qquad&\text{for}\quad i=d+2,\dots,n+2,
\end{cases}
\eeq
where $d=1,\dots,n-1$.
Moreover, the same relations hold for $d=n$, except for \eqref{tiixi} with $i=n+2$,
which is replaced by
\ben
t_{n+2,n+2}(u)\ts \xi_n=\frac{u-n}{u-n-1}\ts\xi_n.
\een
Finally, Proposition~\ref{prop:vectv}\tss(ii) and \eqref{Delta} imply that $\xi_d$
is annihilated by the action of $t_{1j}(u)$ with $j>1$, while
\ben
t_{11}(u)\ts\xi_d=\frac{u+1}{u-d+1}\ts\xi_d.
\een
Thus, the vector $\xi_d$ generates a highest weight module over the algebra
$\X(\osp_{2|2n})$, whose highest weight is found from the action of the series $t_{ii}(u)$.
The formulas for the Drinfeld polynomials easily follow.
\epf

We are now in a position to complete the proof of the Main Theorem.
Proposition~\ref{prop:laon} implies that the
irreducible highest weight
representations of $\X(\osp_{2|2n})$, associated with any tuple
of Drinfeld polynomials of
the first type in \eqref{fundam} is finite-dimensional.
Furthermore, Proposition~\ref{prop:xikdp} shows that the irreducible highest weight
representations of $\X(\osp_{2|2n})$ associated
with the tuples of the form
\eqref{tudp}
and \eqref{tudpn} are also finite-dimensional. On the other hand, by applying
the shift automorphism \eqref{shift} and using the transition rule
\eqref{traru}, we can get all tuples of Drinfeld polynomials
of the second type in \eqref{fundam}
from certain tuples of the first type and those appearing in \eqref{tudp}
and \eqref{tudpn}. This proves the sufficiency of the conditions of the Main Theorem
for the representation $L(\la(u))$ to be finite-dimensional.
The last part of the theorem is an easy consequence of the decomposition
\eqref{tensordecom}; cf. \cite[Cor.~5.19]{amr:rp}.

\appendix

\section{Polynomial evaluation modules over the Yangian $\Y(\gl_{m|n})$}
\label{sec:pe}

\subsection{Defining relations and representations}
\label{subsec:dbA}

Given nonnegative integers $m$ and $n$, we will use the notation $\bi=0$ for $i=1,\dots,m$
and $\bi=1$ for $i=m+1,\dots,m+n$.
Introduce the $\ZZ_2$-graded
vector space $\CC^{m|n}$ over $\CC$ with the basis
$e_1,e_2,\dots,e_{m+n}$, where the parity of the basis vector
$e_i$ is defined to be $\bi\mod 2$.
Accordingly, equip
the endomorphism algebra $\End\CC^{m|n}$ with the $\ZZ_2$-gradation, where
the parity of the matrix unit $e_{ij}$ is found by
$\bi+\bj\mod 2$.

A standard basis of the general linear Lie superalgebra $\gl_{m|n}$ is formed by elements $E_{ij}$
of the parity $\bi+\bj\mod 2$ for $1\leqslant i,j\leqslant m+n$ with the commutation relations
\ben
[E_{ij},E_{kl}]
=\de_{kj}\ts E_{i\tss l}-\de_{i\tss l}\ts E_{kj}(-1)^{(\bi+\bj)(\bk+\bl)}.
\een
The {\em Yang $R$-matrix} associated with $\gl_{m|n}$ is the
rational function in $u$ given by
\beql{ruA}
R(u)=1-P\tss u^{-1},
\eeq
where $P$ is the permutation operator,
\ben
P=\sum_{i,j=1}^{m+n} e_{ij}\ot e_{ji}(-1)^{\bj}\in \End\CC^{m|n}\ot\End\CC^{m|n}.
\een
Following \cite{n:qb} and \cite{n:yg},
define the {\em Yangian}
$\Y(\gl_{m|n})$
as the $\ZZ_2$-graded algebra with generators
$t_{ij}^{(r)}$ of parity $\bi+\bj\mod 2$, where $1\leqslant i,j\leqslant m+n$ and $r=1,2,\dots$,
satisfying the quadratic relations
\ben
\big[\tss t_{ij}(u),t_{kl}(v)\big]=\frac{1}{u-v}
\big(t_{kj}(u)\ts t_{il}(v)-t_{kj}(v)\ts t_{il}(u)\big)
(-1)^{\bi\tss\bj+\bi\tss\bk+\bj\tss\bk},
\een
written in terms of the
formal series
\ben
t_{ij}(u)=\de_{ij}+\sum_{r=1}^{\infty}t_{ij}^{(r)}\ts u^{-r}
\in\Y(\gl_{m|n})[[u^{-1}]].
\een
Combining them into the matrix $T(u)=[t_{ij}(u)]$ and regarding it as the element
\ben
T(u)=\sum_{i,j=1}^{m+n} e_{ij}\ot t_{ij}(u)(-1)^{\bi\tss\bj+\bj}
\in \End\CC^{m|n}\ot \Y(\gl_{m|n})[[u^{-1}]],
\een
we can write the defining relations in the standard $RTT$-form \eqref{RTT}
with the $R$-matrix \eqref{ruA}.

The universal enveloping algebra $\U(\gl_{m|n})$ can be regarded as a subalgebra of
$\Y(\gl_{m|n})$ via the embedding
$
E_{ij}\mapsto t_{ij}^{(1)}(-1)^{\bi},
$
while the mapping
\beql{ev}
t_{ij}(u)\mapsto \de_{ij}+E_{ij}(-1)^{\bi}\ts u^{-1}
\eeq
defines the {\em evaluation homomorphism} $\ev{:}\ts\Y(\gl_{m|n})\to\U(\gl_{m|n})$.

The Yangian $\Y(\gl_{m|n})$ is a Hopf algebra with the coproduct
defined by
\beql{DeltaA}
\De: t_{ij}(u)\mapsto \sum_{k=1}^{m+n} t_{ik}(u)\ot t_{kj}(u).
\eeq

A representation $V$ of the algebra $\Y(\gl_{m|n})$
is called a {\em highest weight representation}
if there exists a nonzero vector
$\xi\in V$ such that $V$ is generated by $\xi$,
\begin{alignat}{2}
t_{ij}(u)\ts\xi&=0 \qquad &&\text{for}
\quad 1\leqslant i<j\leqslant m+n, \qquad \text{and}\non\\
t_{ii}(u)\ts\xi&=\pi_i(u)\ts\xi \qquad &&\text{for}
\quad i=1,\dots,m+n,
\non
\end{alignat}
for some formal series
\beql{laiuA}
\pi_i(u)\in 1+u^{-1}\CC[[u^{-1}]].
\eeq
The vector $\xi$ is called the {\em highest vector}
of $V$ and the $(m+n)$-tuple $\pi(u)=(\pi_{1}(u),\dots,\pi_{m+n}(u))$ is called
its {\em highest weight}.

Given an arbitrary tuple $\pi(u)=(\pi_{1}(u),\dots,\pi_{m+n}(u))$
of formal series of the form \eqref{laiuA},
the {\em Verma module} $M(\pi(u))$ is defined as the quotient of the algebra $\Y(\gl_{m|n})$ by
the left ideal generated by all coefficients of the series $t_{ij}(u)$
with $1\leqslant i<j\leqslant m+n$, and $t_{ii}(u)-\pi_i(u)$ for
$i=1,\dots,m+n$. We will denote by $L(\pi(u))$ its irreducible quotient.
The isomorphism
class of $L(\pi(u))$ is determined by $\pi(u)$.
Necessary and sufficient conditions on $\pi(u)$ for the representation
$L(\pi(u))$ to be finite-dimensional are known due to \cite{zh:sy}.
Their extension to arbitrary parity sequences via odd reflections was
given in \cite{m:or}; see also \cite{l:no}.

Consider the irreducible highest weight representation $V(\pi)$
of the Lie superalgebra $\gl_{m|n}$
with the highest weight $\pi=(\pi_1,\dots,\pi_{m+n})$,
associated with the standard Borel subalgebra.
This means that $E_{ii}\tss\ze=\pi_i\ts\ze$ for $i=1,\dots,m+n$
and $E_{ij}\ts\ze =0$ for $i<j$, where $\ze$ is the highest vector of $V(\pi)$.
The representation $V(\pi)$ is finite-dimensional if and only if the highest weight $\pi$
satisfies the conditions
\ben
\pi_i-\pi_{i+1}\in\ZZ_+
\quad\text{for all}\quad i\ne m;
\een
see \cite{k:rc}.
Use the evaluation homomorphism \eqref{ev} to equip $V(\pi)$ with a $\Y(\gl_{m|n})$-module
structure. The {\em evaluation module} $V(\pi)$
is isomorphic to the Yangian highest weight module $L(\pi(u))$, where
the components of the highest weight $\pi(u)$ are
\ben
\pi_i(u)=1+\pi_i\tss(-1)^{\bi}\ts u^{-1},\qquad i=1,\dots,m+n.
\een

\subsection{Schur--Sergeev duality and fusion procedure}
\label{subsec:sd}

We will follow \cite[Sec.~3.2]{cw:dr} to recall a version of the {\em Schur--Sergeev duality}
going back to \cite{br:hy} and \cite{s:ta}.
An $(m,n)$-{\em hook partition} $\la=(\la_1,\la_2,\dots)$ is a partition with the property
$\la_{m+1}\leqslant n$. This means that the Young diagram $\la$ is contained
in the $(m,n)$-hook as depicted below. The figure also illustrates the partitions
$\mu=(\mu_1,\dots,\mu_m)$ and $\nu=(\nu_1,\dots,\nu_n)$
associated with $\la$. They are introduced by setting
\ben
\mu_i=\max\{\la_i-n,0\},\qquad i=1,\dots,m,
\een
and
\ben
\nu_j=\max\{\la\pr_j-m,0\},\qquad j=1,\dots,n,
\een
where $\la\pr$ denotes the conjugate partition so that $\la\pr_j$ is the length
of column $j$ in the diagram $\la$:

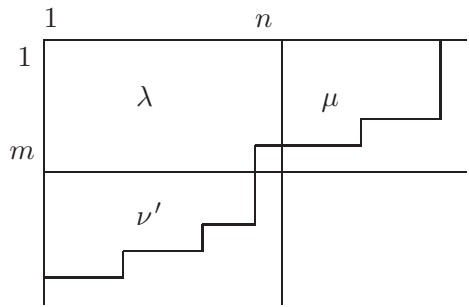
\begin{center}
\begin{picture}(150,110)
\thinlines

\put(0,0){\line(0,1){100}}
\put(90,0){\line(0,1){100}}
\put(0,50){\line(1,0){160}}
\put(0,100){\line(1,0){160}}

\put(0,10){\line(1,0){30}}
\put(30,10){\line(0,1){10}}
\put(30,20){\line(1,0){30}}
\put(60,20){\line(0,1){10}}
\put(60,30){\line(1,0){20}}
\put(80,30){\line(0,1){30}}
\put(80,60){\line(1,0){40}}
\put(120,60){\line(0,1){10}}
\put(120,70){\line(1,0){30}}
\put(150,70){\line(0,1){30}}

\put(0,105){\small$1$}
\put(80,105){\small$n$}
\put(-10,90){\small$1$}
\put(-13,55){\small$m$}

\put(35,75){\small$\la$}
\put(105,75){\small$\mu$}
\put(35,30){\small$\nu\pr$}

\end{picture}
\end{center}

%\vspace{-1cm}

\noindent
We will associate two $(m+n)$-tuples of nonnegative integers with $\la$ by
\ben
\la^{\sharp}=(\la_1,\dots,\la_m,\nu_1,\dots,\nu_n)\Fand
\la^{\flat}=(\mu_1,\dots,\mu_m,\la'_1,\dots,\la'_n).
\een

The tensor product space $(\CC^{m|n})^{\ot d}$ is naturally a module over both $\gl_{m|n}$
and the symmetric group $\Sym_d$. For the action of the basis elements
of $\gl_{m|n}$
we have
\beql{glac}
E_{ij}\mapsto \sum_{a=1}^{d}1^{\ot (a-1)}\ot e_{ij}\ot 1^{\ot (d-a)},
\eeq
while the transposition $(a\ts b)\in\Sym_d$ with $a<b$ acts by $(a\ts b)\mapsto P_{ab}$ with
\ben
P_{ab}=\sum_{i,j=1}^{m+n}1^{\ot (a-1)}\ot e_{ij}\ot 1^{\ot (b-a-1)}
\ot e_{ji}\ot 1^{\ot (d-b)}(-1)^{\bj}.
\een
The images of $\U(\gl_{m|n})$ and $\CC\Sym_d$
in the endomorphism algebra $\End(\CC^{m|n})^{\ot d}$
satisfy the double centralizer property, which leads to the multiplicity-free
decomposition
\ben
(\CC^{m|n})^{\ot d}=\bigoplus_{\la} V(\la^{\sharp})\ot S^{\la},
\een
summed over the $(m,n)$-hook partitions $\la$ with $d$ boxes, where
$S^{\la}$ is the Specht module over $\Sym_d$ associated with $\la$.
By representing the group algebra $\CC\Sym_d$ as the direct sum of matrix algebras
\ben
\CC\Sym_d\cong \underset{\la\vdash d}
\bigoplus\ts\ts\Mat_{f_{\la}}(\CC),\qquad f_{\la}=\dim S^{\la},
\een
we can think of $S^{\la}$ as the canonical
irreducible module over $\Mat_{f_{\la}}(\CC)$ isomorphic to $\CC^{f_{\la}}$.
The diagonal matrix units
$e^{}_{\Uc}=e^{\la}_{\Uc\tss\Uc}\in \Mat_{f_{\la}}(\CC)$ parameterized by
standard $\la$-tableaux $\Uc$ are primitive idempotents in $\CC\Sym_d$.
We may conclude that if $\la$
is an $(m,n)$-hook partition with $d$ boxes, then the image $e^{}_{\Uc}(\CC^{m|n})^{\ot d}$
is isomorphic to the $\gl_{m|n}$-module $V(\la^{\sharp})$. If $\la$ is not contained
in the $(m,n)$-hook, then the image is zero.

Explicit formulas for the idempotents $e^{}_{\Uc}$
can be derived with the use
of the orthonormal Young basis of $S^{\la}$
via the {\em Jucys--Murphy elements\/}
$x_1,\dots,x_d$ of the group algebra $\CC\Sym_d$
defined by
\ben
x_1=0\Fand x_a=(1\ts a)+\dots+(a-1\ts a)\qquad\text{for}\quad a=2,\dots,d.
\een
Given a standard $\la$-tableau $\Uc$, denote by
$\Vc$ the standard tableau
obtained from $\Uc$ by removing the box $\al$ occupied by $d$.
Then the shape of $\Vc$ is a diagram which we denote by $\la^-$.
A box outside $\la^-$ is called {\em addable},
if the union of $\la^-$ and the box is a Young diagram.
We let $c=j-i$ denote the {\em content} of the box $\al=(i,j)$ and let
$a_1,\dots,a_l$ be the contents of all addable boxes of $\la^-$
except for $\al$. The Jucys--Murphy formula gives an inductive
rule for the calculation of $e^{}_\Uc$:
\beql{murphyfo}
e^{}_\Uc=e^{}_\Vc\ts \frac{(x_d-a_1)\dots
(x_d-a_l)}{(c-a_1)\dots (c-a_l)};
\eeq
see \cite{j:fy} and \cite{m:nc}.
The idempotents $e^{}_{\Uc}$ can also be obtained from
the {\em fusion procedure} for the symmetric
group; see \cite{c:ni}, \cite{j:yo} and \cite{n:yc}, and we recall a version following
\cite[Sec.~6.4]{m:yc}, where it was essentially derived from \eqref{murphyfo}.
Take $d$ complex variables $u_1,\dots,u_d$
and consider the rational function with values in $\CC\Sym_d$
defined by
\ben
\phi(u_1,\dots,u_d)=
\prod_{1\leqslant a<b\leqslant d}\Big(1-\frac{(a\ts b)}{u_a-u_b}\Big),
\een
where the product is taken in
the lexicographical order on the set of pairs $(a,b)$.
Suppose that $\la\vdash d$
and let $\Uc$ be a standard $\la$-tableau.
Let $c_a=c_a(\Uc)$ for $a=1,\dots,d$ be the contents of $\Uc$ so that
$c_a=j-i$ if $a$ occupies the box $(i,j)$ in $\Uc$.
Then the consecutive evaluations
of the rational function $\phi(u_1,\dots,u_d)$ are well-defined
and the value coincides with the primitive idempotent $e^{}_\Uc$
multiplied by the product of hook lengths $h(\la)$ of $\la$,
\beql{fus}
\phi(u_1,\dots,u_d)\big|_{u_1=c_1}
\big|_{u_2=c_2}\dots \big|_{u_d=c_d}=h(\la)\ts e^{}_\Uc.
\eeq

\subsection{Yangian action on polynomial modules}
\label{subsec:ya}

The Yang $R$-matrix \eqref{ruA} is a solution of the Yang--Baxter equation
\beql{ybe}
R_{12}(u)\ts R_{13}(u+v)\ts R_{23}(v)=R_{23}(v)\ts R_{13}(u+v)\ts R_{12}(u)
\eeq
in $\End(\CC^{m|n})^{\ot 3}$ with $R_{ab}(u)=1-P_{ab}\ts u^{-1}$.
This implies that
the mapping $T(u)\mapsto R(u)$ defines a representation of the algebra $\Y(\gl_{m|n})$
on the space $\CC^{m|n}$, known as the vector representation. In terms of the generating series it
has the form
\beql{vecreA}
t_{ij}(u)\mapsto \de_{ij}-u^{-1}\tss
e_{ji}(-1)^{\bi\bj}.
\eeq
For an $(m,n)$-hook partition $\la\vdash d$ fix
a standard $\la$-tableau $\Uc$. As above, let $c_1,\dots,c_d$ be the contents
of the respective entries in $\Uc$. By using the coproduct \eqref{DeltaA} and
the shift automorphism of $\Y(\gl_{m|n})$ defined by \eqref{shift}, we get
a representation of the Yangian on the space $(\CC^{m|n})^{\ot d}$ defined by
\beql{tuac}
T(u)\mapsto R_{01}(u-c_1)\dots R_{0d}(u-c_d),
\eeq
which is written in terms of elements of the algebras
\ben
\End\CC^{m|n}\ot\Y(\gl_{m|n})\to \End\CC^{m|n}\ot \End(\CC^{m|n})^{\ot d}
\een
with the first copy of the endomorphism algebra labelled by $0$.

Another form of the vector representation is related to \eqref{vecreA} via
twisting with the super-transposition automorphism
\beql{sutra}
t_{ij}(u)\mapsto t_{ji}(-u)(-1)^{\bi\bj+\bi}.
\eeq
We then get the action of the algebra $\Y(\gl_{m|n})$
on the space $\CC^{m|n}$ given by
\beql{vecreAn}
t_{ij}(u)\mapsto \de_{ij}+u^{-1}\tss
e_{ij}(-1)^{\bi};
\eeq
cf. \eqref{ev}.
It can be written in a matrix form as
$T(u)\mapsto R\pr(-u)$ with
\ben
R\pr(u)=1-Q\ts u^{-1},\qquad Q=
\sum_{i,j=1}^{m+n} e_{ij}\ot e_{ij}\tss(-1)^{\bi+\bj+\bi\bj}.
\een
Accordingly, the composition of the representation \eqref{tuac}
with the automorphism \eqref{sutra} yields another representation
of the Yangian on the space $(\CC^{m|n})^{\ot d}$ given by
\beql{tuactra}
T(u)\mapsto R\pr_{0\tss d}(-u-c_d)\dots R\pr_{01}(-u-c_1),
\eeq
where $R\pr_{0\tss a}(u)=1-Q_{0\tss a}\ts u^{-1}$ and
\ben
Q_{0\tss a}=\sum_{i,j=1}^{m+n} e_{ij}\ot 1^{\ot (a-1)}
\ot e_{ij}\ot 1^{\ot (d-a)}(-1)^{\bi+\bj+\bi\bj}.
\een
Similar to \eqref{tuac}, this representation can also be obtained by using
the opposite coproduct on the Yangian, which is the composition of
\eqref{DeltaA} and the $\ZZ_2$-graded flip operator.

\bth\label{thm:invco}
The subspace $L_{\Uc}=e^{}_{\Uc}(\CC^{m|n})^{\ot d}$
is invariant under the Yangian actions \eqref{tuac} and \eqref{tuactra}.
Moreover, the respective representations of the Yangian $\Y(\gl_{m|n})$ on $L_{\Uc}$
are isomorphic to the highest weight representations $L(\pi^{\flat}(u))$ and $L(\pi^{\sharp}(u))$,
where
\ben
\pi^{\flat}(u)=\big(1-\mu_mu^{-1},\dots,1-\mu_1u^{-1},1+\la'_nu^{-1},\dots,1+\la'_1u^{-1}\big)
\een
and
\ben
\pi^{\sharp}(u)=\big(1+\la_1u^{-1},\dots,1+\la_mu^{-1},1-\nu_1u^{-1},\dots,1-\nu_nu^{-1}\big).
\een
\eth

\bpf
Let $E_{\Uc}\in \End(\CC^{m|n})^{\ot d}$ denote the image of the idempotent $e^{}_{\Uc}$ under
the representation of the symmetric group on $(\CC^{m|n})^{\ot d}$.
We have the relation
\beql{kfus}
R_{01}(u-c_1)\dots R_{0d}(u-c_d)\ts E_{\Uc}
=E_{\Uc}\ts\Big(1-\frac{P_{01}+\dots+P_{0d}}{u}\Big)
\eeq
which plays a key role in the derivation of the fusion formula \eqref{fus};
cf. \cite{n:yc} and \cite[Prop.~6.4.4]{m:yc}.
Apply the anti-automorphism $e_{ij}\mapsto e_{ji}(-1)^{\bi\bj+\bi}$ to the
zeroth copy of the endomorphism algebra $\End\CC^{m|n}$ to derive
\beql{kfustr}
R\pr_{0\tss d}(-u-c_d)\dots R\pr_{01}(-u-c_1)\ts E_{\Uc}
=E_{\Uc}\ts\Big(1+\frac{Q_{01}+\dots+Q_{0d}}{u}\Big).
\eeq
Together with \eqref{kfus} this proves
the first part of the theorem.

Furthermore, relation \eqref{kfustr} shows that the action \eqref{tuactra}
on $L_{\Uc}$ is the composition of the evaluation homomorphism \eqref{ev}
and the action \eqref{glac} of $\gl_{m|n}$. Hence this representation is isomorphic
to the evaluation module $V(\la^{\sharp})\cong L(\pi^{\sharp}(u))$.
Similarly, relation \eqref{kfus}
shows that the Yangian action \eqref{tuac} on the subspace
$L_{\Uc}$ is the composition of the evaluation homomorphism \eqref{ev}
and the action \eqref{glac} of $\gl_{m|n}$ twisted by the automorphism
\beql{psia}
E_{ij}\mapsto -E_{ji}(-1)^{\bi\bj+\bi}.
\eeq
On the other hand, by \cite[Sec.~2.4]{cw:dr}, an application of a chain of
odd reflections
shows that the extremal weight of the $\gl_{m|n}$-module $V(\la^{\sharp})$
with respect to the opposite Borel subalgebra is $\la^{\flat}$. That is, there is a
nonzero vector $\eta\in V(\la^{\sharp})$ of the weight $\la^{\flat}$ such that
$E_{i,i+1}\eta=0$ for all $i\ne m$ and $E_{m+n,1}\eta=0$. By taking the lowest vector
with respect to the action of $\gl_m\oplus\gl_n$ we conclude that $V(\la^{\sharp})$
contains a nonzero vector $\ze$ of the weight $(\mu_m,\dots,\mu_1,\la'_n,\dots,\la'_1)$
such that $E_{ij}\ze=0$ for all $i>j$. Thus, the vector $\ze$ is the highest vector
of the Yangian module $L_{\Uc}$ and its weight is found by taking into account \eqref{ev}
and \eqref{psia} so that this module
is isomorphic to $L(\pi^{\flat}(u))$.
\epf

By using the coproduct \eqref{DeltaA} and the vector representation
\eqref{vecreAn} instead of \eqref{vecreA},
for any complex parameters $z_a$ we get
a representation of the Yangian on the space $(\CC^{m|n})^{\ot d}$ defined by
\ben
T(u)\mapsto R\pr_{01}(-u-z_1)\dots R\pr_{0d}(-u-z_d).
\een
For suitable parameters $z_a$,
the subspaces
$L_{\Uc}=e^{}_{\Uc}(\CC^{m|n})^{\ot d}$
associated with the standard tableaux $\Uc$ of shapes $\la=(d)$ and $\la=(1^d)$,
turn out to be
invariant
under this Yangian action as well.
The primitive idempotents $e^{}_{\Uc}$ associated with the standard
row and column tableaux are, respectively, the symmetrizer
and anti-symmetrizer in $\CC\Sym_d$,
\ben
h^{(d)}=\frac{1}{d\tss !}\sum_{s\in\Sym_d} s
\Fand
a^{(d)}=\frac{1}{d!}\sum_{s\in\Sym_d} \sgn s\cdot s\in\CC\Sym_d.
\een
Note that \eqref{murphyfo} and \eqref{fus} yield multiplicative formulas for
$h^{(d)}$ and $a^{(d)}$.

\bco\label{cor:invcosa}
\quad {\rm (i)}\ \
The subspace $h^{(d)}(\CC^{m|n})^{\ot d}$
is invariant under the action
\ben
T(u)\mapsto R\pr_{01}(-u-d+1)\dots R\pr_{0d}(-u).
\een
Moreover, the representation of the Yangian $\Y(\gl_{m|n})$ on this
subspace is isomorphic to the highest weight representation
with the highest weight
$
\big(1+d\tss u^{-1},1\dots,1\big).
$
\par
{\rm (ii)}\ \
The subspace $a^{(d)}(\CC^{m|n})^{\ot d}$
is invariant under the action
\ben
T(u)\mapsto R\pr_{01}(-u+d-1)\dots R\pr_{0d}(-u).
\een
Moreover, the representation of the Yangian $\Y(\gl_{m|n})$ on this
subspace is isomorphic to the highest weight representation
with the highest weight
\beql{dsm}
\big(\underbrace{1+u^{-1},\dots,1+u^{-1}}_d \ts,1,\dots,1\big)\qquad\text{if}\quad d\leqslant m,
\eeq
and
\ben
\big(\underbrace{1+u^{-1},\dots,1+u^{-1}}_m \ts,1+(m-d)\tss u^{-1},1,\dots,1\big)\qquad\text{if}\quad d> m.
\een
\eco

\bpf
Multiplying both sides
of \eqref{kfustr} by the image $P_{\om}=P_{1, d}\ts P_{2, d-1}\dots$ of the
longest permutation $\om\in\Sym_d$ from the left, we get
\ben
R\pr_{0\tss 1}(-u-c_d)\dots R\pr_{0\tss d}(-u-c_1)\ts P_{\om}\tss E_{\Uc}
=P_{\om}\tss E_{\Uc}\ts\Big(1+\frac{Q_{01}+\dots+Q_{0d}}{u}\Big).
\een
Since $\om\ts h^{(d)}=h^{(d)}$ and $\om\ts a^{(d)}=\sgn\om\cdot a^{(d)}$,
both parts of the corollary follow from the particular cases
of Theorem~\ref{thm:invco} concerning the action \eqref{tuactra}
for the row and column tableaux $\Uc$.
\epf

The following corollary in the case $n=0$ was used in the construction
of the fundamental modules over the Yangians of types $B$, $C$ and $D$ in \cite[Sec.~5.3]{amr:rp}.
It was also applied to the orthosymplectic Yangians in \cite{m:ry} and
in the proof of Proposition~\ref{prop:xikdp} in the previous section.
Equip
the tensor product space $(\CC^{m|n})^{\ot d}$
with the action of $\Y(\gl_{m|n})$ by setting
\ben
t_{ij}(u)\mapsto
\sum_{a_1,\dots,a_{d-1}=1}^{m+n} t_{ia_1}(u-d+1)\ot t_{a_1a_2}(u-d+2)
\ot\dots\ot t_{a_{d-1}j}(u),
\een
where the generators act in the respective copies of the vector space
$\CC^{m|n}$ via the rule \eqref{vecreAn}.
Set
\ben
\xi_d=\sum_{\si\in\Sym_d} \sgn\si\cdot
e_{\si(1)}\ot\dots\ot e_{\si(d)} \in (\CC^{m|n})^{\ot d}.
\een

\bco\label{cor:fungl}
For any $1\leqslant d\leqslant m$ the vector $\xi_d$ has the properties
\ben
t_{ij}(u)\ts \xi_d=0\qquad\text{for}\quad 1\leqslant i<j\leqslant m+n
\een
and
\ben
t_{i\tss i}(u)\ts \xi_d=\begin{cases}
\dfrac{u+1}{u}\ts\xi_d \qquad&\text{for}
\quad i=1,\dots, d,\\[0.5em]
\ \xi_d\qquad&\text{for}\quad i=d+1,\dots,m+n.
\end{cases}
\een
\eco

\bpf
This is immediate from Corollary~\ref{cor:invcosa}\tss(ii), because
the vector $\xi_d$ is the highest vector
of the $\Y(\gl_{m|n})$-module $a^{(d)}(\CC^{m|n})^{\ot d}$
with the highest weight \eqref{dsm}.
\epf

%\newpage
\bigskip\bigskip

\small

\noindent
School of Mathematics and Statistics\newline
University of Sydney,
NSW 2006, Australia\newline
alexander.molev@sydney.edu.au

\end{document}